\documentclass[letterpaper, 10 pt, conference]{ieeeconf}  

\IEEEoverridecommandlockouts                              

\usepackage{graphicx} 
\usepackage{amsmath} 
\usepackage{amssymb}  
\usepackage{subfig}
\usepackage{comment}
\usepackage{cleveref}
\usepackage{algorithm,algorithmic}
\usepackage{xcolor} 
\usepackage[utf8]{inputenc}
\usepackage{balance}
\makeatletter
\newcommand\fs@betterruled{%
  \def\@fs@cfont{\bfseries}\let\@fs@capt\floatc@ruled
  \def\@fs@pre{\vspace*{5pt}\hrule height.8pt depth0pt \kern2pt}%
  \def\@fs@post{\kern2pt\hrule\relax}%
  \def\@fs@mid{\kern2pt\hrule\kern2pt}%
  \let\@fs@iftopcapt\iftrue}
\floatstyle{betterruled}
\restylefloat{algorithm}
\makeatother

\newtheorem{theorem}{Theorem}
\newtheorem{definition}{Definition}

\title{\LARGE \bf
A Consensus-Based Algorithm for Multi-Objective Optimization \\and its Mean-Field Description
}

\author{Giacomo Borghi, Michael Herty and Lorenzo Pareschi
\thanks{This work has been written within the
activities of GNCS group of INdAM (National Institute of
High Mathematics). L.P. acknowledge the partial support of MIUR-PRIN Project 2017, No. 2017KKJP4X “Innovative numerical methods for evolutionary partial differential equations and applications”.  
The work of G.B. is funded by the Deutsche Forschungsgemeinschaft (DFG, German Research Foundation) – Projektnummer 320021702/GRK2326 – Energy, Entropy, and Dissipative Dynamics (EDDy).
}
\thanks{Michael Herty and Giacomo Borghi are with the Institute of Geometry and Practical Mathematics, RWTH Aachen University, 52062 Aachen, Germany
        {\tt\small herty@igpm.rwth-aachen.de}, {\tt\small borghi@eddy.rwth-aachen.de}}%
\thanks{Lorenzo Pareschi and Giacomo Borghi are with the Department of Mathematics and Computer Science, University of Ferrara, 44121 Italy      
	 {\tt\small lorenzo.pareschi@unife.it}}}%

\begin{document}

\maketitle
\thispagestyle{empty}
\pagestyle{empty}

\begin{abstract}
We present a multi-agent algorithm for multi-objective optimization problems, which extends the class of consensus-based optimization methods and relies on a scalarization strategy. The optimization is achieved by a set of interacting agents exploring the search space and attempting to solve all scalar sub-problems simultaneously. We show that those dynamics are described by a mean-field model, which is suitable for a theoretical analysis of the algorithm convergence. Numerical results show the validity of the proposed method.
\end{abstract}



\section{INTRODUCTION}

In applications, decision makers often aim to optimize several objectives which may be in conflict with each others in the sense that improving a solution with respect to an objective may deteriorate another objective. This leads to a so-called multi-objective optimization problem. 

Such problems are often solved by heuristic strategies belonging to the class of evolutionary algorithms \cite{deb2001multi}, where a population of approximations iteratively evolve following mechanisms inspired by natural phenomena. Those are global gradient-free optimization methods which have gained popularity among practitioners, also thanks to their intuitive interpretation as biological systems. Nevertheless, such algorithms typically lack mathematical analysis and understanding \cite{Hooker1995}. For this reason, alternatives are usually developed by
theoretically substantiated methods \cite{Pardalos2018}, leading to a gap between the applications and mathematical research. We refer to \cite{Sergeyev2018} for a discussion on the topic. 

In the single-objective optimization context, a possible bridge between the communities was proposed in \cite{pinnau2017consensus} by developing a Consensus-Based Optimization algorithm (CBO). In \cite{pinnau2017consensus}, the authors place evolutionary algorithms under the more general framework of multi-agent systems, which have received huge attention both in applications \cite{dorri2918multiagents} and in the modeling of biological systems and social interactions \cite{pareschi13}. In CBO methods, several agents interact with each other following a consensus mechanism, collaboratively solving the optimization task. Although the interaction rule is simpler with respect to the most common heuristic strategies, CBO methods are amenable to theoretically analysis as they can be studied by statistical mechanics. First developed to study the dynamics of physical particles in classical mechanics, this mathematical framework allows to give a statistical description of complex systems of agents. Such an approach has been shown to be fruitful to obtain convergence guarantees for single-objective CBO methods \cite{fornasier2021consensusbased} (also for constrained optimization problems \cite{fornasier2020hypersurfaces,carrillo2021constrained,ha2022stochastic,borghi2021constrained}) and for a regularized version of the popular Particle-Swarm Optimization algorithm \cite{grassi2021from, huang2022PSO}. 

In this work, we extend the class of CBO methods to multi-objective optimization problems. The proposed algorithm makes use of a well-known scalarization strategy, the weighted norms approach \cite{book2005mop}, which allows to decompose the multi-objective problem in a set of parametrized single-objective optimization sub-problems. We further couple each agent with a sub-problem and adapt CBO mechanism to solve all sub-problems at the same time. The scalarization strategy and the algorithm are presented in Sections \ref{s:2} and \ref{s:3}. In Section \ref{s:4} we the give a statistical description by presenting the mean-field, approximation of the algorithm dynamics, which is studied analytically. We then test the algorithm on benchmark problems and show the validity of the proposed strategy in Section \ref{s:num}. 

We will use the following notations. With $| \cdot |$ we indicate the $\ell_2$-norm in $\mathbb{R}^n$, and the absolute value on $\mathbb{R}$. Given a Banach space $X$, $\mathcal{P}(X)$ denotes the set of probability measures on $X$, and $\mathcal{P}_q(X)$ the probability measures with bounded $q$-moments.

\section{SCALARIZATION STRATEGY}
\label{s:2}
We are interested in solving a multi-objective optimization problem of the form
\begin{equation}
\min_{x \in \mathbb{R}^d} g(x)  = \left (g_1(x), \cdots, g_m(x) \right)^\top
\label{pb:mop}
\end{equation}

where $g: \mathbb{R}^d \to \mathbb{R}^m$ is a given vector function. We assume $m \geq 2$, as for $m=1$ the problem reduces to a single-objective optimization problem. To define a solution to \eqref{pb:mop}, we rely on the notion of Pareto optimality \cite{jahn2004vector}.

\begin{definition}
A point $\bar x \in \mathbb{R}^d$ is called a \textit{Pareto optimal point} if $g(\bar x)$ is a minimal element of the image set $g(\mathbb{R}^d)$ with respect to the natural partial ordering, that is if there is no $x \in \mathbb{R}^d$ with
\[ g_i(x) \leq g_i(\bar x) \;\; \textup{for all}\;\; i=1, \dots, m\,, \; g(x) \neq g(\bar x)\,.\]
Similarly, $\bar x$ is called a \textit{weakly Pareto optimal point}, if there is no $x \in \mathbb{R}^d$ such that  
\[ g_i(x) < g_i(\bar x)\; \; \textup{for all}\;\; i=1, \dots, m\,.
\]
\end{definition}

In the following, $F_x$ denotes the set of all weakly Pareto points, while $F_g :=g(F_x)$ the so-called weak Pareto front \cite{book2005mop}.

The aim is to find a sufficient number of optimal points to approximate the Pareto front. The proposed method uses a scalarization strategy, which is a common and straightforward way to approach the problem \cite{jahn2004vector,book2005mop}. It allows to solve \eqref{pb:mop} by solving a set of single-objective sub-problems. 

In particular, we consider the approximation sub-problems \cite{jahn2004vector} with the weighted $\ell_p$-norms. Let $\Omega$ be the set of weights vectors
\begin{equation*} \Omega :=\left \{ w \in \mathbb{R}^m_{+}\; | \;\sum_{i=1}^m w_i = 1 \right\}\,.
\end{equation*}
The set $\Omega$ is also known as $m$-dimensional unit, or probability simplex.
Given a $p \in [1,\infty)$, we define the functions $G_p : \mathbb{R}^d \times \Omega \to \mathbb{R}$ as
\begin{equation}
G_p(x,w) := \left( \sum_{k=1}^m w_k |g_k(x)|^p \right)^{\frac1p}
\end{equation}
and extend the definition to the case $p=\infty$
\begin{equation}
G_\infty (x,w) := \max_{k\in \{ 1, \dots, m\}}\, w_k  \,|g_k(x)|\,.
\end{equation}
While being closely related to the $\ell_p$-norms in $\mathbb{R}^m$, we note that  $G_p(\cdot, w)$ is not a norm since, for instance, some weights may be zero. 

Let $p \in [1,\infty]$ be fixed, we then scalarize \eqref{pb:mop}, by considering the sub-problem
\begin{equation}
\min_{x \in \mathbb{R}^d} G_p(x, w)
\label{pb:scal}
\end{equation} 
which is parameterized by a weights vector $w \in \Omega$. The following results hold true.

\begin{theorem}[\cite{jahn2004vector}, Theorem 5.25]
\label{t:1}
Assume $g(x) > 0$ and $p\in[1,\infty]$. For all $w \in \Omega$, any solution to \eqref{pb:scal} is a weakly optimal Pareto point.
\end{theorem}
The inverse implication is not true in general, unless $p=\infty$:
\begin{theorem}[\cite{jahn2004vector}, Corollary 11.21]
Assume $g(x) > 0$ and $p=\infty$. Any weakly optimal Pareto point is a solution \eqref{pb:scal} for a certain $w \in \Omega$ .
\label{t:2}
\end{theorem}

The choice $p=\infty$ is known as the weighted Chebyschev norm approach. By Theorem \ref{t:2}, the Pareto front $F_g$ is approximated by varying the weights vectors in $\Omega$ and solving the correspondent sub-problems \eqref{pb:scal}.  As we will see, the proposed method is derivative-free, so the fact that for $p=\infty$ \eqref{pb:scal} is non-differentiable does not pose an issue. 

\section{M-CBO ALGORITHM}
\label{s:3}

In the following, we propose a CBO algorithm to solve $N \in \mathbb{N}$ sub-problems in the form \eqref{pb:scal}, for an arbitrary (but fixed) choice $p \in [1, \infty]$. Let $\{ w^i\}_{i=1}^N$ be the corresponding weights vectors. To obtain a good approximation of the Pareto front $F_g$, a common choice is to take $\{ w^i\}_{i=1}^N$ uniformly distributed on $\Omega$. We refer to \cite{deb2019generating} for methods to generate uniform points on the unitary simplex given an arbitrary $N$.

Let us start by presenting the classical CBO mechanism and then illustrate how to adapt it to solve all the given sub-problems. 
For this purpose, let $w\in \Omega$ be fixed for the moment and $G(\cdot, w)$ be the correspondent scalar function to minimize.  

As many evolutionary and particle-based methods, CBO algorithms employ a set of $N_p \in \mathbb{N}$ possible solutions and then update their locations repeatedly, according to a specific rule. While we refer to these solutions as agents, sometimes they are also called \textit{particles} see \cite{pinnau2017consensus, carrillo2019consensus}, in view of the parallelism between CBO algorithms and Monte Carlo simulation methods for kinetic equations \cite{pareschi13}. 

Let $\{X^i_k\}_{i=1}^{N_p}$ be 
the agents locations at the $k$-th algorithm iteration. 
In the single-objective CBO method, every agent propagates towards the same point $x_k^\alpha(w)$,
that can be seen as the algorithm approximation to the global minimizer. This value is defined as a convex combination of the agents locations
\begin{equation}
x^\alpha_k(w) := \frac1{Z_\alpha} \sum_{j=1}^{N_p} X^j_k \, \exp\left(-\alpha\, G_p(X^j_k, w)\right)\,,
\label{eq:xa}
\end{equation}
where $\alpha\gg 1$ is a fixed parameter and $Z_\alpha$ is the suitable normalization constant such that the combination is convex.  Taking the limit $\alpha \to \infty$, it holds
\[ x^\alpha_k(w) \;\longrightarrow\;  \underset{i\in\{1,\dots,N\}}{\arg \! \min} G_p(X^i_k, w) \,,
\]
provided that a unique minimum exists. This heuristic strategy promotes the exploration of domain areas where the objective function is lower. The choice of the exponential coefficients  in \eqref{eq:xa} are further justified by the Laplace principle \cite{Dembo2010}, which states that for any absolutely continuous distribution $f \in \mathcal{P}(\mathbb{R}^d)$
\[ \lim_{\alpha \to \infty} - \frac{1}\alpha \log\! \int \!\exp\left (-\alpha G_p(x, w)\right)f(x)dx 
= \hspace{-3mm}\inf_{x \in \textup{supp}(f)}\hspace{-2mm} G(x, w).
\]

After several iterations, the agents concentrate, or \textit{create consensus}, around a point \cite{jin2020convergence}, which is the computed solution to \eqref{pb:scal}.

Instead of iteratively solving $N$ sub-problems, we propose an algorithm which attempts to solve them all simultaneously.  
This can be done assigning a specific sub-problem to every agent, using a specific weights vector $w^i \in \Omega$. To save computation cost, we propose to employ the smallest number of agents to sample all sub-problems, i.e., $N_p=N$. This creates a one-to-one correspondence between the agents and the sub-problems. The method can be easily generalized assigning $n \geq 1$ agents to one sub-problem.

As for the single-objective CBO, a stochastic component is added to the update rule \cite{pinnau2017consensus}.  Two parameters, $\lambda, \sigma >0 $ control the strength of the deterministic and stochastic components, respectively, while $\Delta t>0$ is the step-size.

Let the initial configuration of the agents be independently sampled from a distribution $\rho_0 \in \mathcal{P}(\mathbb{R}^d)$,
\begin{equation}
X_0^i  \sim \rho_0 
\quad \textup{for all} \quad  i=1, \dots,N\,.
\label{eq:init}
\end{equation} 

 The resulting update rule is
\begin{multline}
X^i_{k+1} = X^i_k + \lambda\Delta t\left (x^\alpha_k(w^i) - X^i_k \right)\\
+ \sigma \sqrt{\Delta t} \sum_{l=1}^d (x^\alpha_k(w^i) - X^i_k)_l \, B_{k}^{i,l}\, \vec{e}_l \,,
\label{eq:iter}
\end{multline}
where $(\cdot)_l$ denotes the vector $l$-th component and $\vec{e}_l$ is the unit vector along the dimension $l$. The point $x_k^\alpha(w^i)$ is computed by \eqref{eq:xa} and $B^{i,l}_{k} \in \mathbb{R}$ are randomly sampled values $B^{i,l}_{k} \sim  \mathcal{N}(0, 1)$. Since the update rule is overparametrized, $\lambda$ is usually set to $1$ in CBO algorithms \cite{pinnau2017consensus, carrillo2019consensus}. The Multi-objective CBO method (M-CBO) is summarized in Algorithm \ref{alg}.

In \eqref{eq:iter}, the random component along the $l$-th direction depends on the difference between $(X^i_k)_l$ and $(x_k^\alpha)_l$, such that the exploration behavior is stronger if $i$-the agent is far from the approximation $x^\alpha_k(w^i)$. This type of exploration, hence, is anisotropic \cite{carrillo2019consensus}. Other exploration behaviors proposed for the single-objective CBO algorithm, such as the isotropic exploration \cite{pinnau2017consensus}, can also be considered.

While every point aims to optimize a different sub-problem, all of them are used to compute  $x_k^\alpha$'s, leading to an \textit{interaction} between them. 
Here, the computational cost is saved with the respect to a naive strategy of solving the sub-problems separately. 
Also, heuristically, it exploits the similar structure of the approximation sub-problems.

As a consequence of Theorem \ref{t:1}, we expect the algorithm output $\{X_{\textup{end}}^i\}_{i=1}^N$ to approximate $N$ weakly optimal Pareto points corresponding to the weights vectors:
\begin{equation}\notag
|X_{k}^i - \bar x(w^i)|\xrightarrow{k \to \infty} 0, \quad \textup{with}\quad   g\left (\bar x(w^i) \right) \in F_g\,.
\end{equation}
We verify the convergence in Section \ref{s:num} for two test problems, by measuring the evolution of the average $\ell_2$-error
\begin{equation}
\textup{Err}_{2}(k) := \frac1N \sum_{i=1}^N |X_k^i - \bar x(w^i)|^2\,.
\label{eq:err}
\end{equation} 

Even though $g$ is evaluated only $N$ times per iteration, the overall complexity is $\mathcal{O}(N^2)$ per step. Nevertheless, several random batch techniques developed in the context of particle simulation \cite{AlPa, JLJ} can be used to speed-up the algorithm. The computational complexity is typically reduced to $\mathcal{O}(MN)$ or $\mathcal{O}(M^2)$ for some fixed $M \ll N$. 

\begin{algorithm}
\caption{M-CBO} 
\label{alg}
\begin{algorithmic}
\STATE{Set parameters: $\alpha, \lambda, \sigma, \Delta t$}
\STATE{Initialize the set: $X^i_0 \sim \rho_0\,, i=1, \dots,N$}
\STATE{Select the sub-problems $\{ w^i\}_{i=1}^N $ uniformly in $\Omega$}
\STATE{$k\gets 0$}
\WHILE{stopping criterion is NOT satisfied}
	\STATE{Compute $g(X^i_k)\,, i=1, \dots,N$ }
	\FOR{$i=1, \dots, N$}
		\STATE{compute $x^\alpha(w^i)$ according to \eqref{eq:xa}}
		\STATE{sample $B^{i,l}_{k}$ from $\mathcal{N}(0, 1)$, $l=1, \dots, d$}
		\STATE{compute $X^i_{k+1}$ acording to \eqref{eq:iter} }
	\ENDFOR
	\STATE{$k \gets k+1 $ }
\ENDWHILE
\STATE{$end \gets k $ }
\RETURN $\{ X_{\textup{end}}^i \}_{i=1}^N$
\end{algorithmic}
\end{algorithm}

\section{MEAN-FIELD DESCRIPTION}
\label{s:4}

Mathematically, Algorithm \ref{alg} describes the evolution of $N$ random interacting agents. In this section, we present a mean-field model which approximates such evolution and is amenable to mathematical analysis, following the strategy introduced in \cite{pinnau2017consensus}. The approximation process is twofold. Firstly, the step-size is assumed to be infinitesimal, $\Delta t\to 0$, leading to a continuous-in-time evolution. Secondly, we assume to have infinity agents, $N \to \infty$, leading to a statistical description of the dynamics.

In the context of interacting multi-agent systems, the update rule of Algorithm \ref{alg} typically originates as a simulation of the a continuous-in-time dynamics of $N$ agents \cite{pareschi13}. 
Indeed, \eqref{eq:iter} together with the initial conditions \eqref{eq:init} corresponds to the Euler-Maruyama discretization of the system of stochastic differential equations
\begin{equation}
\begin{split}
dX^i_t &= \lambda \left(x_t^\alpha(w^i) - X^i_t\right)\,dt \\
&\hspace{1.2cm}+  \sigma \sum_{l=1}^d ( x_t^\alpha(w^i) - X^i_t)_l \,dW_{t}^{i,l}\, \vec{e_l}\\
X_0^i &\sim \rho_0 \,,
\end{split}
\label{eq:sde}
\end{equation} 
where $\{(W_t^{i,l})_{t\geq0} \}_{i=1}^N$ are independent Brownian processes. We consider the microscopic model \eqref{eq:sde} to be an approximation of the dynamics \eqref{eq:iter}. 

If $f_0 \in \mathcal P_4(\mathbb{R}^d)$ and $g$ is locally Lipschitz continuous, there exists a unique strong solution $\{ (X^i_t)_{i=1}^N \, | t >0\}$ to \eqref{eq:sde} with continuous almost everywhere paths. 

Since every agent is coupled with a weights vector, 
\begin{equation} 
(X_t^i, w^i) \in \mathbb{R}^d \times \Omega\quad \textup{for all} \quad  i=1, \dots,N\,,
\label{eq:couples}
\end{equation}
the dynamics \eqref{eq:sde} can also be considered to take place in the augmented space $\mathbb{R}^d \times \Omega$. Instead of studying the trajectories of all $N$ different couples $(X_t^i, w^i)$, it is convenient to give a statistical description of the ensemble evolution through a probability density $f$ where -- informally -- $f(t,x,w)$ quantifies the probability of finding at the time $t$ an agent with weights $w \in \Omega$ at the location $x \in \mathbb{R}^d$. 

The empirical distribution $f^N$ associated to the couples \eqref{eq:couples} is
\begin{equation}
f^N(t, x,w) =\frac 1 N \sum_{i=1}^N \delta({X_t^i}-x) \, \delta({w^i}-w)\, ,
\label{eq:emp}
\end{equation}
where $\delta$ is the Dirac distribution. For all $t\geq0$, $f^N(t, \cdot, \cdot)$ is a random variable in $\mathcal{P}(\mathbb{R}^d \times \Omega)$. 

For large values of $N$, a standard strategy is to approximate $f^N$ by its limit as $N \to \infty$, the mean-field limit \cite{bolley2011stochastic}. 

Since the agents are initially sampled from $\rho_0 \in \mathcal{P}(\mathbb{R}^d)$ and the weights vectors are uniformly distributed on $\Omega$, it holds
\begin{equation}
f_0^N \rightharpoonup \rho_0 \otimes \mu \quad \textup{in law as} \quad N\ \to \infty\,,
\end{equation}
where  $f_0^N = f^N(t, \cdot, \cdot)$ and $\mu$ denotes the uniform distribution over $\Omega$. 

The mean-field limit $f$ of the empirical distribution solves the non-linear Fokker-Planck equation
\begin{equation}
\begin{split}
\partial_t f(t,x,w) &= \nabla_x \cdot \left( (x_t^\alpha(w) - x) f(t, x, w) \right) \\
&\hspace{-5mm}+ \frac{\sigma^2}2 \sum_{l=1}^d \frac{\partial^2}{\partial x^2_l} \left(( x^\alpha_t(w) - x)_l^2 f(t, x, w)  \right)\\
\lim_{t\to 0} f(t,x,w) &= \rho_0(x) \mu(w)\,.
\label{eq:mf}
\end{split}
\end{equation}

The system \eqref{eq:mf} describes the mean-field model, which approximates the algorithm computation for a large number $N$ of agents. 
The points $x^\alpha_t (w)$ are now defined as the mean-field equivalent of \eqref{eq:xa}:
\begin{equation}
x_t^\alpha(w) =\frac1Z_\alpha \int_{\mathbb{R}^d} x\, \exp(-\alpha G_p(x,w))\, \rho(t,x)\, dx\,,
\label{eq:xa2}
\end{equation} 
where $\rho(t, x) = \int_{\Omega} f(t,x,w)\, dw$ is the first marginal of $f$.
We note that, even if we initialize the algorithm assigning more then one agent per sub-problem, the correspondent mean-field model is still given by \eqref{eq:mf}.

%

As for Algorithm \ref{alg}, we expect the couple $(x,w)$ to converge towards $(\bar x(w),w)$, with $ \bar x(w) \in F_x$ or, equivalently, that 
\begin{equation}\notag
f(T, x,w) \approx \bar f(x,w), \quad \bar f(x,w):= \delta( \bar x(w) - x) \mu(w)
\end{equation}
at some time $T>0$.  By definition, the first marginal $\bar \rho$ of $\bar f$ is then concentrated on weakly optimal Pareto points. By Theorem \ref{t:2} and $p=\infty$, the weakly optimal points fulfill
\begin{equation*}
\textup{supp} \left(\bar \rho \right) = F_x\,.
\end{equation*}

Whether $f$ actually converges numerically towards $\bar f$ can be analytically studied by the time evolution of 
\begin{equation}
\textup{Err}^{M\!F}_2 (t) = \int_{\mathbb{R}^d}\int_{\Omega} |x-\bar x(w)|^2\, f(t,x,w)\,dw\,dx\,.
\end{equation}
Indeed, the above error is an upper bound for the 2-Wasserstein distance between $f$ and $\bar f$, which metrizes the weak convergence on $\mathcal{P}_2(\mathbb{R}^d \times \Omega)$ \cite{santambrogio2015optimal}. 
Such an approach has been shown successful to prove the convergence of single-objective CBO methods \cite{fornasier2021consensusbased}. We leave this type of analysis of the proposed multi-objective optimization method for future work. 

\section{NUMERICAL RESULTS}
\label{s:num}

\begin{figure}[!b]
\centering
\includegraphics[trim=6cm 10cm 6cm 10.2cm, clip, width=3in]{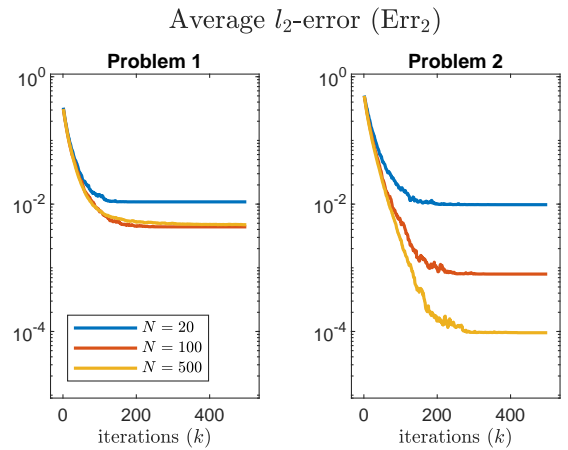}
\caption{Average $\ell_2$-error, as defined in \eqref{eq:err}, as function of the iterative step $k$ for Problem 1 and 2. Different population sizes $N$ are considered. Update rule given by \eqref{eq:iter}. Results averaged over 1000 runs.}
\label{fig:conv}
\end{figure}

In this section, we validate the suggested heuristic strategy on four bi-objective optimization test problems, that is with $m=2$. 
The first problem \cite{Pardalos2018} is given by 
\begin{align}
g(x_1,x_2) = \begin{pmatrix}
5(x_1 - 0.1)^2 + (x_2 - 0.1)^2
\\
(x_1 - 0.9)^2 + 5(x_2 - 0.9)^2
\end{pmatrix}
\end{align}
with $x \in [0,1]^2$, while the second one is the DEB2DK problem \cite{finding2004deb}, with one knee ($K=1$) and $d=2$. Problems 3 and 4 are taken form the test suit \cite[problems UF4, UF7]{zhang2008multiobjective} and can be scaled to different dimensions $d$. Given that Problems 2 and 3 are non-convex, we employ the Chebyschev norm approach by setting $p = \infty$ which, thanks to Theorem \ref{t:2}, is suitable for such problems. 

As in many multi-objective problems \cite{Pardalos2018}, the search space of each test consists of a different hypercube $D \subset \mathbb{R}^d$. Therefore, after every iteration we clip the vectors $X_k^i$ to ensure the agents will remain in the search space. Alternatively, a penalization term can be added to the objective function. The initial agents positions are uniformly sampled from $D$, $\rho_0 = \mathcal{U}(D)$, while the weights are deterministically chosen as $w^i = \left((i-1)\Delta w, 1 - (i-1)\Delta w \right)$ with $\Delta w = 1/(N-1)$.



\begin{figure*}[!t]
\centering
\subfloat[$d=2,\, N=100,\, \sigma = 4$]{\includegraphics[trim=7cm 10cm 7.3cm 10cm, clip, width=1.6in]{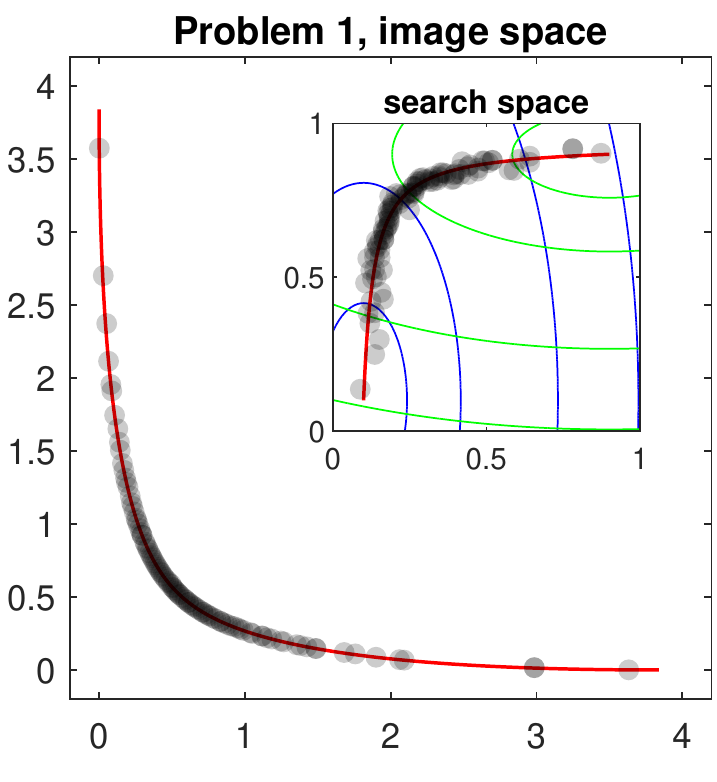}
\label{fig2a}
}
\hfil
\subfloat[$d=2,\, N=100,\, \sigma = 4$]{\includegraphics[trim=7cm 10cm 7.3cm 10cm, clip, width=1.6in]{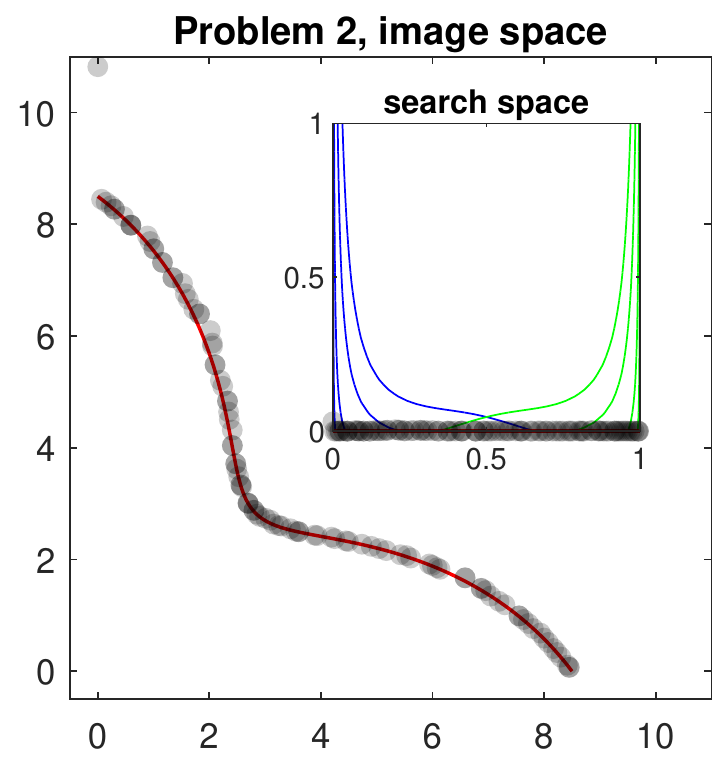} \label{fig2b}}
\hfil
\subfloat[$d=5,\, N=300,\, \sigma = 10$]{\includegraphics[trim=7cm 10cm 7.3cm 10cm, clip, width=1.6in]{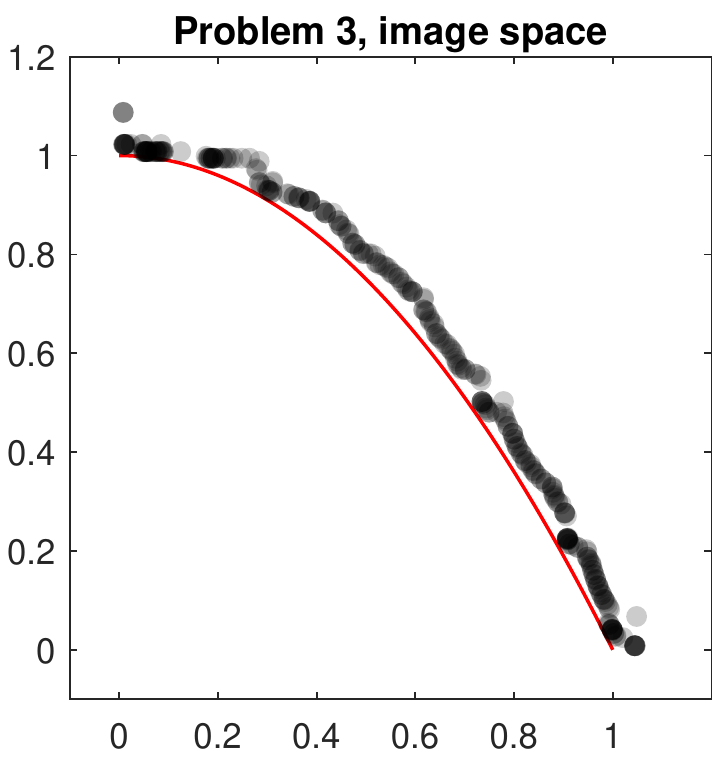}\label{fig2c}} 
\hfil
\subfloat[$d=5,\, N=300,\, \sigma = 10$]{\includegraphics[trim=7cm 10cm 7.3cm 10cm, clip, width=1.6in]{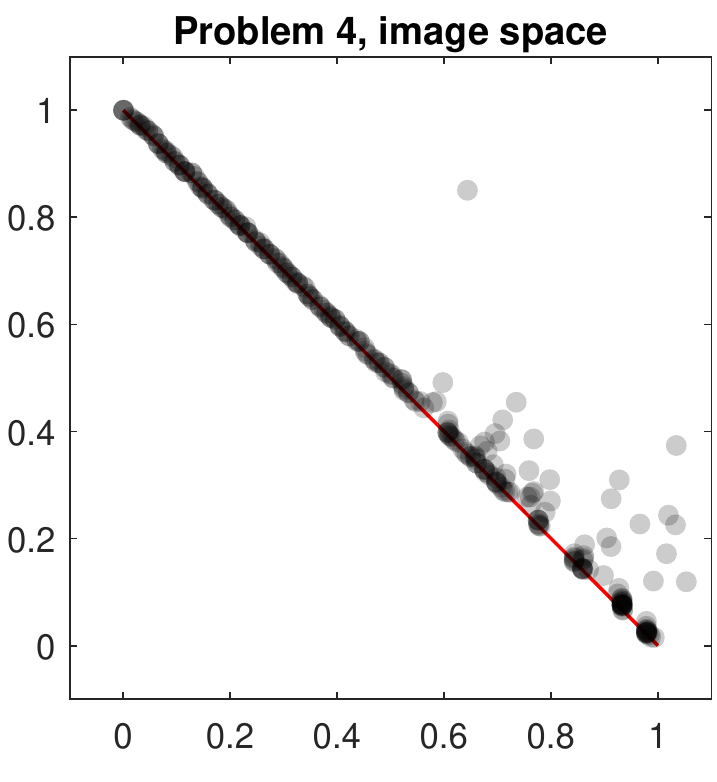}\label{fig2d}} 
\caption{Agents position at the end of a single run, $k_{\max}=500$ . In red the Pareto front $F_g$ and, when $d=2$, the optimal points $F_x$. In Problems 3 and 4 the greedy strategy \eqref{eq:greedy} is used.}
\label{fig:pareto}
\end{figure*}

In similar multi-agent algorithms, the step size $\Delta t$ is typically of order between $0.1$ and $0.01$ \cite{pinnau2017consensus, carrillo2019consensus, benfenati2021binary}. In our experiments, we were able to reach higher accuracy with step size $\Delta t =0.01$.
Even though using adapting parameters is a common strategy to promote exploration at the beginning of the computation \cite{carrillo2019consensus,fhps20-2}, we test the algorithm mechanism by keeping $\lambda =1, \Delta t =0.01, \alpha = 10^{5}$ fixed. The exploration parameter is fixed to $\sigma=4$ for Problems 1 and 2, and $\sigma=10$ for higher dimensional problems. The agents evolve until a maximum number of iterations $k_{\textup{max}} $ is reached.

We first consider Problem 1 and 2, for which is possible to compute the solution of each sub-problem \eqref{pb:scal} with high accuracy.
Figs. \ref{fig:conv} show the average error $\textup{Err}_2$ for different population sizes $N$. We note that the agents converge quickly to the solution of the correspondent sub-problems, up to a certain accuracy. 

For the remaining part of the computation, the agent system is almost stationary. Indeed, once the $i$-th agent approximates the solution of the $i$-th problem better than the other agents, the convex combination $x_k^\alpha(w^i)$ is typically close to $X_k^i$ itself due to the Laplace principle. Note that, even though using small values of $\alpha$ would still evolve the dynamics, this would lead to a movement of positions towards the same point and an increase of the average error. This due the definition of $x_k^\alpha(w)$ as convex combination of the agents position. 


As expected, employing larger number of agents allows to a better exploration of the search space and higher accuracy of the solution. Nevertheless, this improvement may be little, as shown in Fig. \ref{fig:conv} (left), where the accuracy reached with $N=100$ agents is almost as good as the one  reached with $N=500$. This can be explained by a known shortcoming of the scalarization strategy. Typically, taking uniformly distributed weights vectors does not guarantee equally distributed Pareto points, both in the search-space and in image-space \cite{Pardalos2018}.

This is particularly evident for the Problem 1, see Fig. \ref{fig2a}, where the solutions are concentrated at the center (with respect to the symmetry axis). Indeed, the agents concentrate towards the center due to the combinatorial nature of $x^\alpha(w)$. This leads to a low accuracy 
solutions at the extrema of the Pareto front,
even for large population sizes. In Problem 2, this effect is less evident as the agents are better distributed over the front, see Fig. \ref{fig2b}. 

\begin{figure}[!b]
\centering
\includegraphics[trim=6cm 10cm 6cm 10.2cm, clip, width=3in]{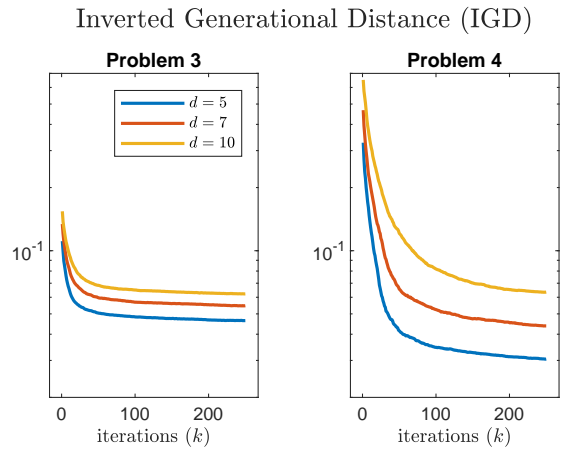}
\caption{IGD metric as function of the iterative step $k$ for Problem 3 and 4. Different dimensions $d=5,7,10$ of the search space are considered, while $N=300$ is fixed.  Update rule given by \eqref{eq:greedy} (greedy strategy). Results averaged over 50 runs.}
\label{fig:dim}
\end{figure}

We test the proposed method with higher dimensional problems, where $d=5,7,10$. To improve the  algorithm performance, we employ a greedy strategy which updates the location of the $i$-th agent only if the new location improves the objective function value of the $i$-th sub-problem. This can be done by substituting the update rule \eqref{eq:iter} with
\begin{align}
Y^i_{k+1} &= X^i_k + \lambda\Delta t\left (x^\alpha_k(w^i) - X^i_k \right) \notag
\\
&\hspace{1cm} + \sigma \sqrt{\Delta t} \sum_{l=1}^d (x^\alpha_k(w^i) - X^i_k)_l \, B_{k}^{i,l}\, \vec{e}_l 
\label{eq:greedy}
\\
X^i_{k+1} & =H\left( G_p(X^i_k, w^i) - G_p(Y^i_{k+1}, w^i)  \right) ( Y^i_{k+1} - X^i_k)  \notag \\
&\hspace{1cm} + X^i_k \notag
\end{align}
where $H$ is the Heaviside function ($H(x) = 0$ if $x\leq0$ and $1$ otherwise). We note that the above dynamics can also be described by a correspondent mean-field model, after an appropriate regularization of $H$. 
For simplicity, we omit it and refer to \cite{pinnau2017consensus, fornasier2021consensusbased} for more details. We remark the method might benefit from a mixed strategy where the greedy mechanism is gradually turned on, as in the Simulated Annealing algorithm \cite{kirkpatrick1983annealing}. We leave this investigation for future work.

Fig. \ref{fig:dim} shows the evolution of the Inverted Generational Distance (IGD) which is a common metric for multi-objective optimization tasks, as it measures both the convergence and the well-distribution of the points over the Pareto front \cite{cheng2012performance}. Even though the IGD metric is increasing together with the dimension, the method is able to converge towards the Pareto front with only $N=300$ agents in dimension $d=10$.

\section{CONCLUSIONS}
We proposed a multi-agent multi-objective algorithm which solves, at the same time, $N$ sub-problems generated by a scalarization strategy. We give a statistical description of the agents dynamics through the mean-field limit, which is given by the limit of the step-size $\Delta t \to 0$ and $N \to \infty$. Such approximation is an essential step to analytically study the algorithm's behavior.
Computational tests show the agents distribution over the Pareto front and the validity of the proposed approach.

The initial choice of uniform weights vectors is known to be non-optimal for a general multi-objective problem.  We observe that this may have a strong impact on the convergence. In future work, we plan of adding an interaction in the weights to mitigate this issue and obtain a better approximation of the Pareto front and, also, an higher accuracy in solution of the sub-problems.

\IEEEtriggeratref{13}

\bibliographystyle{IEEEtran}
\bibliography{IEEEabrv,IEEEexample}

\end{document}